\documentclass[12pt]{amsart}

\usepackage[english]{babel}
\usepackage[utf8]{inputenc}
\usepackage{times}
\usepackage[T1]{fontenc}
\usepackage{graphicx}
\usepackage{amscd,amsmath,amsfonts,amstext,amssymb,amsbsy,amsopn,amsthm,eucal}
\usepackage{txfonts}
\usepackage[figuresright]{rotating}
\usepackage{longtable,booktabs}
\usepackage{rotating,booktabs}
\usepackage{caption}
\usepackage{multirow}
\usepackage{tensor}
\usepackage[table]{xcolor}
\usepackage{fancyhdr}
\usepackage{color} 
\usepackage{hyperref}

\hypersetup{
	pdftitle={Talenti-Robin},
	pdfauthor={CW},
	pdfsubject={Eigenvalue},
	pdfkeywords={Bossel-Daners, Isoperimetric inequality, Robin Laplacian, First eigenvalue}
	pdfpagelayout=SinglePage,
	pdfpagemode=UseOutlines,  pdfpagemode=UseOutlines,
	colorlinks,
	linkcolor=[rgb]{0,0,0.7},
	urlcolor=[rgb]{0,0,0.4},
	citecolor=[rgb]{0.4,0.1,0}
}

\newcommand{\R}{\mathbb{R}^n}
\newcommand{\Sp}{\mathbb{S}^n}

\newcommand{\abs}[1]{\left|#1\right|}

\newcommand{\ak}{\alpha_\kappa}
\newcommand{\akz}{\alpha_0}

\newcommand{\n}{\nabla_g}

\newcommand{\ld}{\lambda_{1,\beta}(\Omega)}
\newcommand{\ldb}{\lambda_{1,\beta}(\Omega^\sharp)}
\newcommand{\snk}{\operatorname{sn_\kappa}}

\DeclareMathOperator{\Ric}{Ric}

\def\dfrac{\displaystyle\frac}

\newtheorem{prop}{Proposition}[section]
\newtheorem{thm}[prop]{Theorem}

\newtheorem{lem}[prop]{Lemma}
\newtheorem{rem}[prop]{Remark}
\newtheorem{cor}[prop]{Corollary}
\newtheorem{defn}[prop]{Definition}

\numberwithin{equation}{section}

\begin{document}
	
	
\baselineskip=17pt
	
	
\title[Talenti comparison Robin]{ Comparison results for solutions of Poisson equations with Robin boundary on complete Riemannian manifolds}

\author{Daguang Chen,Haizhong Li, Yilun Wei}

\thanks{The authors were supported by NSFC-FWO 11961131001 and  NSFC grant No. 11831005.}
	
\subjclass[2010]{{53C44}, {53C42}}
\keywords{Bossel-Daners inequality, Isoperimetric inequality, Talenti's comparison, Lorentz norm}

\maketitle

	

\begin{abstract}	
	 In this paper, by using Schwarz rearrangement and isoperimetric inequalities, we prove comparison results for the solutions of Poisson equations on complete Riemannian manifolds with  $Ric\geq (n-1)\kappa$, $\, \kappa\geq 0$, which extends the results in \cite{ANT-Talenti-a}. Furthermore, as applications of our comparison results, we obtain the Saint-Venant inequality and Bossel-Daners inequality for Robin Laplacian.  
\end{abstract}
	
\section{Introduction}
Let $(M,g)$ be an $n$-dimensional complete Riemannian manifold with $\Ric(g)\geq (n-1)\kappa$, where $\kappa\geq 0$. For $\kappa=0$, we further assume that $M$ is noncompact with positive asymptotic volume ratio, i.e.  
$$
\operatorname{AVR}(g)=\lim\limits_{r\rightarrow\infty}\frac{|B_r(x_0)|_{g}}{\omega_n r^n }>0,
$$ 
where $B_r(x_0)$ denotes the geodesic ball in $M$ centered at $x_0$ with radius $r$, while $\omega_n$ denotes the volume of $n$-dimensional Euclidean unit ball. Let $\Omega$ be a bounded domain with smooth boundary $\partial\Omega$ on $M$.  For a positive number $\beta$ and a nonnegative
function $f$ in $L^2(\Omega)$, we consider the following problem
\begin{equation}\label{eq1}
	\left\{
	\begin{array}{ll}
		-\Delta_{g} u= f,  & \text{in}\quad  \Omega,\\\\
		\dfrac{\partial u}{\partial N} +\beta u =0, & \text{on}\quad \partial\Omega,
	\end{array}
	\right.
\end{equation}
where $N$ denotes the outer unit normal to $\partial\Omega$. For $\beta=\infty$, the equation \eqref{eq1} can be seen as  the problem equipped  with  Dirichlet boundary condition.

According to isoperimetric inequalities in \cite{Brendle} for $\kappa=0$ (See also \cite{AFM,BK202012}) and in \cite{Gromov} for $\kappa=1$, we define   
\begin{equation}\label{iso-cont}
	\ak=
	\begin{cases}
		\operatorname{AVR}(g), & \kappa=0,\\
		\frac{\abs{M}_{g}}{\abs{\mathbb{S}^n}}, & \kappa=1,
	\end{cases}
\end{equation}
where $\abs{M}_g$, $\abs{\Sp}$ denote the volume of $M$ and $\Sp$, respectively. Let $(\R(\kappa),g_\kappa)$ be Euclidean space when $\kappa=0$ and be unit sphere when $\kappa=1$ with canonical metrics. In this article, we intend to establish a comparison principle with the solution of the following problem
\begin{equation}\label{eq2}
	\left\{
	\begin{array}{ll}
		-\Delta_{g_\kappa} v=f^\sharp,  & \text{in}\quad \Omega^\sharp,\\\\
		\frac{\partial v}{\partial N}+\beta v=0, & \text{on}\quad \partial\Omega^\sharp,
	\end{array}
	\right.
\end{equation}
where $\Omega^\sharp$ is a geodesic ball on $\R(\kappa)$ satisfying $|\Omega|_{g}=\ak|\Omega^\sharp|_{g_\kappa}$ and $f^\sharp$ is the Schwarz rearrangement of $f$.

 For Poisson equations with Dirichlet boundary conditions in Euclidean spaces,  Talenti \cite{Talenti76} gave pointwise comparisons of $u^\sharp$ and $v$. Talenti's comparison results were generalized to semilinear and nonlinear elliptic equations, for instance, in \cite{Talenti79, ALT90, ColladayLM18,Mondin0Vedovato}. Recently,  Talenti's comparison results were extended to  solutions of Poisson equations on complete noncompact Riemannian manifolds with nonnegative Ricci curvature by  the first two authors in \cite{CL2021}.  We also refer the reader to excellent books \cite{Baernstein19,Kawohl,Kesavan2006} for related topics.

 Recently, Alvino et al. \cite{ANT-Talenti-a} obtained the Talenti type comparison results of Poisson equations with Robin boundary conditions in Euclidean spaces.   The results in \cite{ANT-Talenti-a} were  generalized to the case of $p$-Laplace 
 operator with Robin boundary conditions in \cite{AGM} and to the torsion problem 
 for the Hermite operator with Robin boundary 
 conditions in \cite{CGNT2021}. Following the strategies in \cite{ANT-Talenti-a}, we focus on Poisson equations with Robin boundary conditions on complete Riemannian manifolds  and prove the following comparison results.
\begin{thm}[]\label{thm:main1}
	Let $(M,g)$ be an $n$-dimensional complete Riemannian manifold with $\Ric(g)\geq (n-1)\kappa$, where $\kappa=0$ or $1$. For $\kappa=0$, we further assume that $M$ is noncompact with $\operatorname{AVR}(g)>0$. Let $\Omega$ be a bounded domain with smooth boundary on $M$ and $\Omega^\sharp$ be a geodesic ball on $\R(\kappa)$ satisfying $|\Omega|_{g}=\ak|\Omega^\sharp|_{g_\kappa}$, where $\ak$ is defined in \eqref{iso-cont}. Let $u$ and $v$ be solutions to \eqref{eq1} and \eqref{eq2} respectively, then
\begin{equation}\label{thm:main1-1}
		||u||_{L^{p,1}(\Omega)}\leq \ak^{\frac{1}{p}}||v||_{L^{p,1}(\Omega^\sharp)}\quad \text{for} ~ \kappa=0,1 ~ \text{and} ~ 0< p\leq\frac{n}{2n-2},
	\end{equation}
\begin{equation}\label{thm:main1-2}
		||u||_{L^{2p,2}(\Omega)}\leq \ak^{\frac{1}{2p}}||v||_{L^{2p,2}(\Omega^\sharp)}\quad \text{for} ~ \kappa=0 ~ \text{and} ~  0< p\leq \frac{n}{3n-4},
	\end{equation}
\begin{equation}\label{thm:main1-4}
	||u||_{L^{2p,2}(\Omega)}\leq \ak^{\frac{1}{2p}}||v||_{L^{2p,2}(\Omega^\sharp)}\quad \text{for} ~ \kappa=1 ~ \text{and} ~ 0< p\leq \frac{n}{3n-3}.
\end{equation}
In particular, when $n=2$, $\kappa=1$ and $0< p\leq 1$,
\begin{equation}\label{thm:main1-5}
	||u||_{L^{2p,2}(\Omega)}\leq \ak^{\frac{1}{2p}}||v||_{L^{2p,2}(\Omega^\sharp)}.
\end{equation}
\end{thm}

\begin{thm}[]\label{thm:main2}
Under the same assumptions as Theorem \ref{thm:main1}, letting  $u$ and $v$ be solutions to \eqref{eq1} and \eqref{eq2} for $f\equiv 1$, respectively, we have
\begin{equation}\label{thm:main2-1}
	||u||_{L^{p,1}(\Omega)}\leq \ak^{\frac{1}{p}}||v||_{L^{p,1}(\Omega^\sharp)}\quad \text{for} ~ \kappa=0,1 ~ \text{and} ~ p\leq\frac{n}{n-2},
\end{equation}
and
\begin{equation}\label{thm:main2-2}
	||u||_{L^{2p,2}(\Omega)}\leq \ak^{\frac{1}{2p}}||v||_{L^{2p,2}(\Omega^\sharp)}\quad \text{for} ~ \kappa=0 ~ \text{and} ~ 0< p\leq\frac{n}{n-2}.
\end{equation}
In particular, for $n=2$ and  $\kappa=0$ , we have a pointwise comparison result
\begin{equation}\label{est-point}
	u^\sharp(x)\leq v(x) \quad \text{for all} ~ x\in \Omega^\sharp.
\end{equation}
\end{thm}

Let $(M, g)$ be a complete Riemannian manifold and let $\Omega\subset M$ be a smoothly bounded  domain in $M$.  Assume that $u$ is a positive solution to the following  torsion problem with Robin boundary condition
\begin{equation}\label{equ:tor}
	\left\{
	\begin{array}{ll}
		-\Delta_g u= 1,  & \text{in $\Omega$},\\\\
		\dfrac{\partial u}{\partial N} +\beta u =0, & \text{on $\partial\Omega$}.
	\end{array}
	\right.
\end{equation}
The  torsional rigidity $T_\beta(\Omega)$ with Robin boundary condition  of the domain $\Omega$ is defined by
\[
T_\beta(\Omega) = \int_{\Omega} u\, dV_g.
\]

In 2015, by using free discontinuity techniques, Bucur et al. in \cite{BG15} proved the Saint-Venant inequality for Robin Laplacian in Euclidean spaces. In 2019,  Alvino et al. in \cite{ANT-Talenti-a} obtained the same Saint-Venant inequality for Robin Laplacian via a Talenti type comparison result in Euclidean spaces. 

Choosing $p=1$ in  \eqref{thm:main2-1}, we  can deduce the Saint-Venant inequality for Robin Laplacian on given manifolds.
\begin{cor}[]\label{cor:tor}
Under the same assumptions as Theorem \ref{thm:main1}, for the torsion problem \eqref{equ:tor},  we obtain
	\begin{equation}\label{cor:tor-1}
		T_\beta(\Omega) \leq \ak T_\beta(\Omega^\sharp).
	\end{equation}
\end{cor}
\begin{rem}
	For the torsion problem with Dirichlet boundary condition on  Riemannian manifold $(M,g)$ satisfying $\Ric(g)\geq (n-1)$, the inequality \eqref{cor:tor-1} is due to \cite{ColladayLM18} and \cite{GHM2015}. In \cite{CL2021}, the first two  authors obtained  sharp estimates for the $L^1$-moment spectrum and  $L^\infty$-moment spectrum with Dirichlet boundary conditions on bounded domains in complete Riemannian manifolds with nonnegative Ricci curvature.
\end{rem}


 Let $\Omega$ be a bounded domain with smooth boundary on $M$ and $\Omega^\sharp$ be a geodesic ball on $\R(\kappa)$ satisfying $|\Omega|_{g}=\ak|\Omega^\sharp|_{g_\kappa}$. Let $\ld$ and $\ldb$ be the first eigenvalues to \eqref{eig-1} and \eqref{eig-2} respectively, i.e.
\begin{equation}\label{eig-1}
	\left\{
	\begin{array}{ll}
		-\Delta_{g} u= \ld u,  & \text{in $\Omega$},\\\\
		\dfrac{\partial u}{\partial N} +\beta u =0, & \text{on $\partial\Omega$},
	\end{array}
	\right.
\end{equation}
and
\begin{equation}\label{eig-2}
	\left\{
	\begin{array}{ll}
		-\Delta_{g_\kappa} v=\ldb v,  & \text{in $\Omega^\sharp$},\\\\
		\dfrac{\partial v}{\partial N}+\beta v=0, & \text{on $\partial\Omega^\sharp$}.
	\end{array}
	\right.
\end{equation}
As another application to Theorem \ref{thm:main1}, we use Talenti comparisons to deal with Bossel-Daners inequality for $n=2$. When $n>2$, however, $p=1$ no longer satisfies the conditions of \eqref{thm:main1-2} and \eqref{thm:main1-5}, hence we fail to use Theorem \ref{thm:main1} to give a direct proof. 
\begin{thm}[]\label{thm:eig}
	Let $(M,g)$ be an $n$-dimensional complete Riemannian manifold with $\Ric(g)\geq (n-1)\kappa$, where $\kappa=0$ or $1$. For $\kappa=0$, we further assume that $M$ is noncompact with $\operatorname{AVR}(g)>0$. Let $\Omega$ be a bounded domain with smooth boundary on $M$ and $\Omega^\sharp$ be geodesic ball on $\R(\kappa)$ satisfying $|\Omega|_{g}=\operatorname{\alpha_\kappa}|\Omega^\sharp|_{g_\kappa}$. Then
	\begin{equation}\label{thm:eig-1}
	\lambda_{1,\beta}(\Omega)\geq \lambda_{1,\beta}(\Omega^\sharp).
\end{equation}	
Moreover, the inequality holds in \eqref{thm:eig-1} if and only if $(M,g)$ is isometric to $(\R(\kappa),g_\kappa)$ and $\Omega$ is isometric to a geodesic ball in $(\R(\kappa),g_\kappa)$.
\end{thm}

\begin{rem}The inequality \eqref{thm:eig-1} is also called Faber-Krahn inequality  for Robin Laplacian.  When $(M,g)$ is Euclidean space, the result is due to Bossel in \cite{Bossel} for $n=2$ and Daners in \cite{D1} for all dimensions (See also \cite{DaiFu2011,BD,BG15}). For bounded domains in Riemannian manifolds with $\Ric(g)\geq n-1$, the inequality \eqref{thm:eig-1} is due to Chen-Cheng-Li in \cite{CCL2021}.
\end{rem}

The paper is organized as following. In Section \ref{Pre}, we recall the isoperimetric inequality of complete Riemannian manifolds with  $Ric\geq (n-1)\kappa$, where $\kappa\geq 0$,  Schwarz rearrangement, Lorentz space and Gronwall's inequality; 
In Section \ref{lemmas}, we establish some integral inequalities for the solutions to \eqref{eq1} and \eqref{eq2} by using the isoperimetric inequality;
In Section \ref{Proof of MT}, we give the proofs of Theorem \ref{thm:main1} and Theorem  \ref{thm:main2}; In Section \ref{pf:BD}, we give the proof of Theorem \ref{thm:eig}.

\section{Preliminaries}\label{Pre}

\subsection{Isoperimetric inequalities}

Let $(M,g)$ be an $n$-dimensional complete Riemannian manifold with $\Ric(g)\geq (n-1)\kappa$, where $\kappa=0$ or 1. For $\kappa=0$, we further assume that $M$ is non-compact with $\operatorname{AVR}(g)>0$. Let $\Omega$ be a bounded domain with smooth boundary on $M$ and $\Omega^\sharp$ be a geodesic ball on $\R(\kappa)$ satisfying $|\Omega|_{g}=\ak|\Omega^\sharp|_{g_\kappa}$. According to Brendle  \cite{Brendle} (See also \cite{AFM,BK202012}) and L\'{e}vy-Gromov \cite{Gromov}, we have the isoperimetric inequality
\begin{equation}\label{isoperi}
	\abs{\partial \Omega}_{g}\geq \ak \abs{\partial \Omega^\sharp}_{g_\kappa}.
\end{equation}
The equality  holds in \eqref{isoperi} if 
and only if $(M,g)$ is isometric to $(\R(\kappa),g_\kappa)$ and $\Omega$ is isometric to a geodesic ball in $\R(\kappa)$.  

For a complete noncompact $n$-dimensional Riemannian manifold $(M,g)$  with  nonnegative Ricci curvature and positive asymptotic volume growth, we also notice that the inequality \eqref{isoperi}  is proved by Agostiniani et al. in \cite[Theorem 1.8]{AFM} for $n=3$ and then extended  to $3\leq n\leq 7$ by Fogagnolo and Mazzieri in \cite{FM202012}. Furthermore, \eqref{isoperi} and its equality case still hold in ${\sf CD} (0,N)$  metric measure spaces based on the method of optimal mass transport by Balogh and Krist\'{a}ly in \cite{BK202012}.

\subsection{Schwarz rearrangement}
Let $u$ and $v$ be solutions to \eqref{eq1} and \eqref{eq2} respectively.  For $t\geq 0$ we denote by 
\begin{equation*}
	\begin{aligned}
	&U_t=\{x \in \Omega: u(x)>t\}, \quad \partial U_t^i = \partial U_t \cap \Omega, \quad \partial U_t^e= \partial U_t \cap \partial\Omega, \quad \mu(t) = |U_t|_{g},\\
	&V_t=\{x \in \Omega^\sharp: v(x)>t\}, \quad
	\phi(t) = |V_t|_{g_\kappa}.	
	\end{aligned}	
\end{equation*}

Denoting by $u_m$ and $v_m$ the minimum of $u$ and $v$ respectively, thanks to the positiveness of $\beta$ and Robin boundary conditions, we have $u_m \geq 0$ and $v_m \geq 0$. Since $v$ is radial, positive and decreasing along the radius, we have $V_t$ coincides with $\Omega^\sharp$ when $0\leq t<v_m$.
\begin{defn}
	Letting $h: \Omega \to \mathbb{R}$ be a measurable function, the distribution function of $h$ is the function $\mu_h : [0,+\infty)\, \to [0, +\infty)$ defined by
	$$
	\mu_h(t)= |\{x \in \Omega \, :\,  |h(x)|>t\}|_{g}.
	$$
	\end{defn}
\begin{defn}	
	The decreasing rearrangement $h^*:[0,|\Omega|_{g}]\to \mathbb{R}$ is defined based on the distribution function of $h$, that is
	\begin{equation}
		h^*(s)=
		\begin{cases}
			\underset{\Omega}{\textup{esssup}}\,h, & s=0,\\
			\inf \{t:\mu_h(t)\leq s \}, & s>0.
		\end{cases}
	\end{equation}
\end{defn}
\begin{defn}\label{SchwarzS}	
	The Schwarz rearrangement $h^\sharp:\Omega^\sharp\to \mathbb{R}$ is defined based on the decreasing rearrangement of $h$, that is
	\begin{equation}
		h^\sharp(x)= h^*(\ak|B^\kappa_R|_{g_\kappa}),
	\end{equation}
where $B^\kappa_R$ is the geodesic ball in $\R(\kappa)$ with radius $R$.
	\end{defn}
The distribution function of $h$ and $h^\sharp$ satisfies
\begin{equation}\label{Eqn:RearrMean}
	\mu_{h}(t)=\ak\mu_{h^{\sharp}}(t).
\end{equation}
By definition, $h$, $h^*$ and $h^\sharp$ are equi-distributed in the sense that
\begin{equation}\label{equi}
	\displaystyle{||h||_{L^p(\Omega)}=||h^*||_{L^p(0, |\Omega|_{g})}=\ak^{\frac{1}{p}}||h^\sharp||_{L^p(\Omega^\sharp)}}.
\end{equation}
Given measurable functions $h_1,h_2$ on $\Omega$, the Hardy-Littlewood inequality holds,
\begin{equation}\label{inq:HL}
	\int_{\Omega}{|h_1(x)h_2(x)| dV_{g}} \le \int_{0}^{|\Omega|_{g}}{h_1^*(s) h_2^*(s) ds}.
\end{equation}
Choosing $h_2=\chi_{\left\lbrace |u|>t\right\rbrace}$ in \eqref{inq:HL}, one has
\begin{equation}\label{inq:re}
	\int_{|u|>t}{|h_1(x)| dV_{g}} \le \int_{0}^{\mu(t)}{h_1^*(s) ds}.
\end{equation}

By strong maximum principle, both solutions $u$ and $v$  to \eqref{eq1} and \eqref{eq2}  achieve their minimum on boundaries. Hence $u$ and $v$ are strictly positive in the interior of domains. Moreover, by isoperimetric inequality \eqref{isoperi}, we have
\begin{equation}
	\begin{split}
		\ak v_m  |\partial\Omega^\sharp|_{g_\kappa} &= \ak\int_{\partial \Omega^\sharp} v(x) \, d\mu_{g_\kappa}= \frac{\ak}{\beta}\int_{\Omega^\sharp} f^\sharp \, dV_{g_\kappa}=\frac{1}{\beta} \int_{\Omega} f \, dV_{g} \\
		& = \int_{\partial \Omega} u(x) \, d\mu_{g} \\
		&\geq u_m  |\partial\Omega|_{g} \\
		&\geq u_m \ak |\partial\Omega^\sharp|_{g_\kappa},
	\end{split}
\end{equation}
which implies that 
\begin{equation}\label{minima_eq}
	 \min_\Omega u=u_m \leq v_m=\min_{\Omega^\sharp} v .
\end{equation}
An important consequence of \eqref{minima_eq} is that
\begin{equation} \label{mf}
	\mu (t) \leq |\Omega|_{g}=\ak \phi (t) \quad \text{for} \quad 0\leq t< v_m.
\end{equation} 

\subsection{Lorentz space}

\begin{defn}\label{defn:Lorentz}
Let $0<p<+\infty$ and $0<q\le +\infty$. The Lorentz space $L^{p,q}(\Omega)$ is the space of functions such that the quantity
	\begin{equation*}
		||h||_{L^{p,q}} =
		\begin{cases}
			\displaystyle{ p^{\frac{1}{q}} \left( \int_{0}^{\infty}  t^{q} \mu_h(t)^{\frac{q}{p}}\ \frac{dt}{t}\right)^{\frac{1}{q}}}, & \text{if}~~ 0<q<\infty,\\
			\displaystyle{\sup_{t>0} \, (t^p \mu_h(t))}, & \text{if}~~ q=\infty,
		\end{cases}
	\end{equation*}
is finite.
\end{defn}

It's well known that Lorentz space coincides with $L^p$ space when $p=q$ (See \cite{Talenti94} for more details). 

\subsection{Gronwall's inequality}
Let $\xi$ be a continuously differentiable function satisfying $$\tau\xi'(\tau)\leq \xi(\tau)+C$$ for all $\tau \geq\tau_0>0$, where $C$ is a non-negative constant. For all $\tau\geq \tau_0$, we have
\begin{equation}\label{Gronwall's inequality}
	\xi'(\tau)\leq\frac{\xi(\tau_0)+C}{\tau_0}.
\end{equation}

\section{Some lemmas}\label{lemmas}
Let $(M,g)$ be an $n$-dimensional complete Riemannian manifold with $\Ric(g)\geq (n-1)\kappa$, where $\kappa=0$ or 1. Moreover, for $\kappa=0$, assume $M$ is non-compact with $\operatorname{AVR}(g)>0$. Let $\Omega$ be a bounded domain with smooth boundary on $M$ and $\Omega^\sharp$ be a geodesic ball on $\R(\kappa)$ satisfying $|\Omega|_{g}=\ak|\Omega^\sharp|_{g_\kappa}$. Define a function $G_\kappa(l)$,
\begin{equation}\label{Gk}
	G_\kappa(l)=\frac{dI_\kappa}{dr}\circ I_\kappa^{-1}(l),
\end{equation}
where 
$$I_\kappa(r)=n\omega_n\ak \int_{0}^{r}{\snk^{n-1}(s) ds},$$
and	
\begin{equation*}
	\snk(s)=
	\begin{cases}
		s, & \kappa=0,\\
		\sin s, & \kappa=1.
	\end{cases}
\end{equation*}
We then have the following lemma:

\begin{lem}
	
Let u and v be solutions to \eqref{eq1} and \eqref{eq2} respectively. For almost every $t>0$, we have
\begin{equation}\label{lem1-1}
	G_\kappa(\mu(t))^2\leq \int_0^{\mu(t)}{f^*(s) ds}\cdot \left(-\mu'(t)+\frac{1}{\beta}\int_{\partial U_t^e}{\frac{1}{u}d\mu_{g}}\right),
\end{equation}
	
\begin{equation}\label{lem1-2}
	\tilde{G}_\kappa(\phi(t))^2= \ak^{-1}\int_0^{\ak\phi(t)}{f^*(s) ds}\cdot\left(-\phi'(t)+\frac{1}{\beta}\int_{\partial V_t\cap\partial\Omega^\sharp}{\frac{1}{v}d\mu_{g_\kappa}}\right),
\end{equation}
where $\tilde{G}_\kappa(r)=\ak^{-1}G_\kappa(\ak r).$
\end{lem}

\begin{proof}
		Multiplying \eqref{eq1} by $\varphi\in H^1(\Omega)$ and integrating by parts, we have 
\begin{equation}\label{10}
	\begin{aligned}
		\int_{\Omega}{\nabla_{g} u\cdot\nabla_{g} \varphi dV_{g}}+\beta\int_{\partial\Omega}{u\varphi d\mu_{g}}
		&=\int_{\Omega}{\nabla_{g} u\cdot\nabla_{g} \varphi dV_{g}}+\int_{\partial\Omega}{(-\varphi) \frac{\partial u}{\partial N}d\mu_{g}}\\
		&=\int_{\Omega}{\left[\nabla_{g} u\cdot\nabla_{g} \varphi-div_{g}(\varphi\nabla_{g} u) \right]dV_{g}}\\
		&=\int_{\Omega}{-\varphi \Delta_{g} u dV_{g}}\\
		&=\int_{\Omega}{f\varphi dV_{g}}.
	\end{aligned}
\end{equation}
Define a test function, for $h>0$,
\begin{equation*}
	\varphi_h=
	\begin{cases}
		0, & 0<u\leq t,\\
		h,& u>t+h,\\ 
		u-t,& t<u\leq t+h.\\
	\end{cases}
\end{equation*}
Choosing $\varphi=\varphi_h$ in \eqref{10}, one has
\begin{equation}
	\begin{aligned}
		&\int_{U_t\backslash U_{t+h}}{|\nabla_{g} u|_{g}^2 dV_{g}}+\beta\int_{\partial U_{t+h}^e}{uh d\mu_{g}}+\beta\int_{\partial U_{t}^e\backslash \partial U_{t+h}^e}{u(u-t) d\mu_{g}}\\
		&=\int_{U_t\backslash U_{t+h}}{f(u-t) dV_{g}}+\int_{U_{t+h}}{fh dV_{g}}.
	\end{aligned}
\end{equation}
Dividing by $h$ and letting $h\rightarrow 0+$, we have
\begin{equation}
	-\frac{d}{dt}\left(\int_{U_t}{|\nabla_{g} u|_{g}^2} dV_{g}\right)+\beta\int_{\partial U_t^e}{u d\mu_{g}}=\int_{U_t}{f dV_{g}}.
\end{equation} 
Applying co-area formula, we then obtain
\begin{equation}\label{13}
	\int_{\partial U_t^i}{|\nabla_{g} u|_{g} d\mu_{g}}+\int_{\partial U_t^e}{\beta u d\mu_{g}}=\int_{U_t}{f dV_{g}}.
\end{equation}
Define $$A=\left\{
	\begin{array}
		{l@{\quad\quad}l}
		|\nabla_{g} u|_{g}, & \partial U_t^i, \\
		\beta u, & \partial U_t^e. \\
	\end{array}
	\right.$$
Then \eqref{13} becomes
\begin{equation}\label{14}
	\int_{\partial U_t}{A d\mu_{g}}=\int_{U_t}{f dV_{g}}.
\end{equation} 
By \eqref{isoperi}, \eqref{inq:re}, \eqref{Gk}, \eqref{14} and co-area formula, we have
\begin{equation}\label{3.10}
	\begin{aligned}
		G_\kappa(\mu(t))^2&\leq|\partial U_t|_{g}^2\\
		&\leq\int_{\partial U_t}{A d\mu_{g}}\cdot\int_{\partial U_t}{A^{-1} d\mu_{g}}\\
		&=\int_{U_t}{f dV_{g}}\cdot\left(\int_{\partial U_t^i}{\frac{1}{|\nabla_{g} u|_{g}} d\mu_{g}}+\int_{\partial U_t^e}{\frac{1}{\beta u} d\mu_{g}}\right)\\
		&=\int_{U_t}{f dV_{g}}\cdot\left(-\mu'(t)+\frac{1}{\beta}\int_{\partial U_t^e}{\frac{1}{u} d\mu_{g}}\right)\\
		&\leq\int_0^{\mu(t)}{f^*(s) ds}\cdot\left(-\mu'(t)+\frac{1}{\beta}\int_{\partial U_t^e}{\frac{1}{u} d\mu_{g}}\right).\\		
	\end{aligned}
\end{equation}
Similarly, by $\tilde{G}_\kappa(\phi(t))=|\partial V_t|_{g_\kappa}$ and $\int_{V_t}{f^\sharp dV_{g_\kappa}}=\ak^{-1}\int_{0}^{\ak\phi(t)}{f^*(s) ds}$, we prove \eqref{lem1-2}.
	\end{proof}
 
\begin{lem}
	Let u and v be solutions to \eqref{eq1} and \eqref{eq2} respectively.  For all $t>v_m$, we have
	\begin{equation}\label{lem2-1}
		\int_{0}^{t}{\tau\left(\int_{\partial U_\tau^e}{\frac{1}{u}d\mu_{g}}\right)d\tau}\leq \frac{1}{2\beta}\int_{0}^{|\Omega|_{g}}{f^*(s) ds},
	\end{equation}
\begin{equation}\label{lem2-2}
	\int_{0}^{t}{\tau\left(\int_{\partial V_\tau\cap\partial\Omega^\sharp}{\frac{1}{v}d\mu_{g_\kappa}}\right)d\tau}= \frac{\ak^{-1}}{2\beta}\int_{0}^{|\Omega|_{g}}{f^*(s) ds}.
\end{equation}
	\end{lem} 
\begin{proof}
	By Fubini's theorem, \eqref{eq1} and \eqref{equi}, we have 
\begin{equation}
	\begin{aligned}
		\int_{0}^{+\infty}{\tau\left(\int_{\partial U_\tau^e}{\frac{1}{u}d\mu_{g}}\right)d\tau}&=\int_{0}^{+\infty}{d\tau\left(\int_{\partial \Omega\cap\{u\geq\tau\}}{\frac{\tau}{u}d\mu_{g}}\right)}\\
		&=\int_{0}^{+\infty}{d\tau\left(\int_{\partial \Omega}{\chi_{\{u\geq\tau\}}\frac{\tau}{u}d\mu_{g}}\right)}\\
		&=\int_{\partial \Omega}{\left(\int_{0}^{u}{\tau d\tau}\right)\frac{1}{u}d\mu_{g}}\\
		&=\frac{1}{2}\int_{\partial \Omega}{u d\mu_{g}}\\
		&=\frac{1}{2\beta}\int_{\Omega}{f dV_{g}}\\
		&=\frac{1}{2\beta}\int_{0}^{|\Omega|_{g}}{f^*(s) ds}.
	\end{aligned}
\end{equation}
Thus, for all $t> v_m$, we obtain
	\begin{equation}
	\begin{aligned}
		\int_{0}^{t}{\tau\left(\int_{\partial U_\tau^e}{\frac{1}{u}d\mu_{g}}\right)d\tau}&\leq \int_{0}^{+\infty}{\tau\left(\int_{\partial U_\tau^e}{\frac{1}{u}d\mu_{g}}\right)d\tau}\\
		&=\frac{1}{2\beta}\int_{0}^{|\Omega|_{g}}{f^*(s) ds}.
	\end{aligned}
	\end{equation}
With the fact that$\int_{\Omega^\sharp}{f^\sharp dV_{g_\kappa}}=\ak^{-1}\int_{0}^{|\Omega|_{g}}{f^*(s) ds}$ and $\partial V_t\cap \partial\Omega^\sharp=\emptyset$ when $t> v_m$, we deduce
\begin{equation}
	\begin{aligned}
		\int_{0}^{t}{\tau\left(\int_{\partial V_\tau\cap \partial\Omega^\sharp}{\frac{1}{v}d\mu_{g_\kappa}}\right)d\tau}&= \int_{0}^{+\infty}{\tau\left(\int_{\partial V_\tau\cap \partial\Omega^\sharp}{\frac{1}{v}d\mu_{g_\kappa}}\right)d\tau}\\
		&=\frac{\ak^{-1}}{2\beta}\int_{0}^{|\Omega|_{g}}{f^*(s) ds}.
	\end{aligned}
\end{equation}
\end{proof}

\section{Proofs of Theorem \ref{thm:main1} and Theorem  \ref{thm:main2}}\label{Proof of MT}
In this section, we give the proof of main Theorem \ref{thm:main1} and Theorem  \ref{thm:main2} by using the lemmas in previous section.
\begin{proof}[Proof of Theorem \ref{thm:main1}]
		Multiplying \eqref{lem1-1} by $t\mu(t)^{\frac{1}{p}}G_\kappa(\mu(t))^{-2}$ and integrating from $0$ to $\tau>v_m$, we have
\begin{equation}\label{52}
	\begin{aligned}
		\int_{0}^{\tau}{t\mu(t)^{\frac{1}{p}} dt}\leq& \int_{0}^{\tau}{-t\mu'(t) \mu(t)^{\frac{1}{p}} G_\kappa(\mu(t))^{-2} \int_{0}^{\mu(t)}{f^*(s) ds} dt}\\
		&+\int_{0}^{\tau}{\left[\frac{t}{\beta}\int_{\partial U_t^e}{\frac{1}{u}d\mu_{g}}\right] \mu(t)^{\frac{1}{p}} G_\kappa(\mu(t))^{-2} \int_{0}^{\mu(t)}{f^*(s) ds} dt}.
	\end{aligned}
\end{equation}
From \eqref{Gk}, it's obvious that $l^{\frac{1}{p}}G_0(l)^{-2}=n^{-2}(\omega_n\operatorname{\alpha_0})^{-\frac{2}{n}}l^{\frac{1}{p}-2+\frac{2}{n}}$ is non-decreasing when $0<p\leq \frac{n}{2n-2}$. Besides, $l^{\frac{1}{p}}G_1(l)^{-2}$ is also non-decreasing when $0<p\leq \frac{n}{2n-2}$. By definition of $G_1$ in \eqref{Gk}, we have
\begin{equation}
	\begin{aligned}
		\frac{d}{dl}\left(l^{\frac{1}{p}}G_1(l)^{-2}\right)&=\frac{1}{p}l^{\frac{1}{p}-1}\cdot \frac{G_1(l)-2plG_1'(l)}{G_1(l)^3}\\
		&=\frac{p^{-1}l^{\frac{1}{p}-1}}{(I_1'\circ I_1^{-1}(l))^4}\cdot (I_1'^2-2pI_1I_1'')\circ I_1^{-1}(l).
	\end{aligned}
\end{equation}
It is sufficient to show $k(r):=I_1'(r)^2-2pI_1(r)I_1''(r)\geq 0$ on $[0, \pi)$. Noticing that $$\cos r\cdot\int_{0}^{r}{\sin^{n-1}s\, ds}\cdot\sin^{-n}r\leq \frac{1}{n}$$ on $[0, \pi)$, one has
\begin{equation}
	\begin{aligned}
		k(r)&=\left(n\omega_n\operatorname{\alpha_1} \sin^{n-1} r\right)^2\left[1-2p(n-1)\frac{\cos r\cdot\int_{0}^{r}{\sin^{n-1}s\, ds}}{\sin^nr}\right]\geq 0		
	\end{aligned}
\end{equation}
when $0<p\leq \frac{n}{2n-2}$. By monotonicity of
$l^{\frac{1}{p}}G_\kappa(l)^{-2}$ when $0<p\leq \frac{n}{2n-2}$ and \eqref{lem2-1}, it follows
\begin{equation}\label{4.2}
	\begin{aligned}
		&\int_{0}^{\tau}{\left[\frac{t}{\beta}\int_{\partial U_t^e}{\frac{1}{u}d\mu_{g}}\right] \mu(t)^{\frac{1}{p}} G_\kappa(\mu(t))^{-2} \int_{0}^{\mu(t)}{f^*(s) ds} dt}\\
		&\leq |\Omega|_{g}^{\frac{1}{p}} G_\kappa(|\Omega|_{g})^{-2} \int_{0}^{|\Omega|_{g}}{f^*(s) ds}\cdot \int_{0}^{\tau}{\frac{t}{\beta}\int_{\partial U_t^e}{\frac{1}{u}d\mu_{g}} dt}\\
		&\leq \frac{|\Omega|_{g}^{\frac{1}{p}} G_\kappa(|\Omega|_{g})^{-2}}{2\beta^2}\left(\int_{0}^{|\Omega|_{g}}{f^*(s) ds}\right)^2.
	\end{aligned}
\end{equation}
By \eqref{52} and \eqref{4.2}, we have
\begin{equation}\label{54}
	\begin{aligned}
		\int_{0}^{\tau}{t\mu(t)^{\frac{1}{p}} dt}\leq& -\int_{0}^{\tau}{t\left(\mu'(t) \mu(t)^{\frac{1}{p}} G_\kappa(\mu(t))^{-2} \int_{0}^{\mu(t)}{f^*(s) ds} \right) dt}\\
		&+\frac{|\Omega_{g}|^{\frac{1}{p}} G_\kappa(|\Omega|_{g})^{-2}}{2\beta^2}\left(\int_{0}^{|\Omega|_{g}}{f^*(s) ds}\right)^2\\
		=&-\int_{0}^{\tau}{tdF_\kappa(\mu(t))}+\frac{|\Omega|_{g}^{\frac{1}{p}} G_\kappa(|\Omega|_{g})^{-2}}{2\beta^2}\left(\int_{0}^{|\Omega|_{g}}{f^*(s) ds}\right)^2,
	\end{aligned}
\end{equation}	
where $F_\kappa(l)=\int_{0}^{l}{w^{\frac{1}{p}}G_\kappa(w)^{-2}\int_{0}^{w}{f^*(s) ds} dw}$.
Integrating by parts and applying \eqref{Gronwall's inequality} , we obtain
\begin{equation}\label{55}
	\begin{aligned}
		\int_{0}^{\tau}{\mu(t)^{\frac{1}{p}} dt}+F_\kappa(\mu(\tau))\leq\frac{1}{v_m}&\Bigg[\int_{0}^{v_m}{dt \int_{0}^{t}{\mu(r)^{\frac{1}{p}} dr} }+\int_{0}^{v_m}{F_\kappa(\mu(t)) dt}\\&+\frac{|\Omega|_{g}^{\frac{1}{p}} G_\kappa(|\Omega|_{g})^{-2}}{2\beta^2}\left(\int_{0}^{|\Omega|_{g}}{f^*(s) ds}\right)^2\Bigg].
	\end{aligned}
\end{equation}
Making similar computations and using \eqref{lem1-2} and \eqref{lem2-2}, we have
\begin{equation}\label{56}
	\begin{aligned}
		\int_{0}^{\tau}{t\phi(t)^{\frac{1}{p}} dt}=& -\int_{0}^{\tau}{t\left(\operatorname{\alpha_\kappa}^{-1}\phi'(t) \phi(t)^{\frac{1}{p}} \tilde{G}_\kappa(\phi(t))^{-2} \int_{0}^{\operatorname{\alpha_\kappa}\phi(t)}{f^*(s) ds} \right) dt}\\
		&+\operatorname{\alpha_\kappa}^{-1}\int_{0}^{\tau}{\left[\frac{t}{\beta}\int_{\partial V_t\cap \partial \Omega^\sharp}{\frac{1}{v}d\mu_{g_\kappa}}\right] \phi(t)^{\frac{1}{p}} \tilde{G}_\kappa(\phi(t))^{-2} \int_{0}^{\operatorname{\alpha_\kappa}\phi(t)}{f^*(s) ds} dt}\\
		=&-\int_{0}^{\tau}{td\left(\operatorname{\alpha_\kappa}^{-\frac{1}{p}}F_\kappa(\operatorname{\alpha_\kappa} \phi(t))\right)}+\operatorname{\alpha_\kappa}^{-\frac{1}{p}}\frac{|\Omega|_{g}^{\frac{1}{p}} G_\kappa(|\Omega|_{g})^{-2}}{2\beta^2}\left(\int_{0}^{|\Omega|_{g}}{f^*(s) ds}\right)^2.
	\end{aligned}	
\end{equation}
Thus,
\begin{equation}\label{57}
	\begin{aligned}
		\int_{0}^{\tau}{(\operatorname{\alpha_\kappa}\phi(t))^{\frac{1}{p}} dt}+F_\kappa(\operatorname{\alpha_\kappa}\phi(\tau))=\frac{1}{v_m}&\Bigg[\int_{0}^{v_m}{dt \int_{0}^{t}{(\operatorname{\alpha_\kappa}\phi(r))^{\frac{1}{p}} dr} }+\int_{0}^{v_m}{F_\kappa(\operatorname{\alpha_\kappa}\phi(t)) dt}\\
		&+\frac{|\Omega|_{g}^{\frac{1}{p}} G_\kappa(|\Omega|_{g})^{-2}}{2\beta^2}\left(\int_{0}^{|\Omega|_{g}}{f^*(s) ds}\right)^2\Bigg].
	\end{aligned}
\end{equation}
Since \eqref{mf} holds, by \eqref{55} and \eqref{57}, we conclude
\begin{equation}
	\int_{0}^{\tau}{\mu(t)^{\frac{1}{p}} dt}+F_\kappa(\mu(\tau))\leq \int_{0}^{\tau}{(\operatorname{\alpha_\kappa}\phi(t))^{\frac{1}{p}} dt}+F_\kappa(\operatorname{\alpha_\kappa}\phi(\tau)).
\end{equation}
Letting $\tau\rightarrow +\infty$, we prove \eqref{thm:main1-1}.

Next, we come to prove \eqref{thm:main1-2}, \eqref{thm:main1-4} and \eqref{thm:main1-5}. 
Letting $\tau\rightarrow +\infty$ in \eqref{54} and \eqref{56}, we have
\begin{equation}\label{4.10}
		\int_{0}^{+\infty}{t\mu(t)^{\frac{1}{p}} dt}\leq\int_{0}^{+\infty}{F_\kappa(\mu(t)) dt}+\frac{|\Omega|_{g}^{\frac{1}{p}} G_\kappa(|\Omega|_{g})^{-2}}{2\beta^2}\left(\int_{0}^{|\Omega|_{g}}{f^*(s) ds}\right)^2,
\end{equation}
\begin{equation}\label{4.11}
	\int_{0}^{+\infty}{t(\operatorname{\alpha_\kappa}\phi(t))^{\frac{1}{p}} dt}=\int_{0}^{+\infty}{F_\kappa(\operatorname{\alpha_\kappa}\phi(t)) dt}+\frac{|\Omega|_{g}^{\frac{1}{p}} G_\kappa(|\Omega|_{g})^{-2}}{2\beta^2}\left(\int_{0}^{|\Omega|_{g}}{f^*(s) ds}\right)^2.
\end{equation}
We then need to prove 
$$\int_{0}^{+\infty}{F_\kappa(\mu(t)) dt}\leq \int_{0}^{+\infty}{F_\kappa(\operatorname{\alpha_\kappa}\phi(t)) dt}.$$ Multiplying \eqref{lem1-1} by $tF_\kappa(\mu(t))G_\kappa(\mu(t))^{-2}$ and integrating from $0$ to $\tau>v_m$, we get
\begin{equation}
	\begin{aligned}
		\int_{0}^{\tau}{tF_\kappa(\mu(t)) dt}\leq& -\int_{0}^{\tau}{t\left(\mu'(t) F_\kappa(\mu(t)) G_\kappa(\mu(t))^{-2} \int_{0}^{\mu(t)}{f^*(s) ds} \right) dt}\\
		&+\int_{0}^{\tau}{\left[\frac{t}{\beta}\int_{\partial U_t^e}{\frac{1}{u}d\mu_{g}}\right] F_\kappa(\mu(t)) G_\kappa(\mu(t))^{-2} \int_{0}^{\mu(t)}{f^*(s) ds} dt}.
	\end{aligned}
\end{equation}	
It's not hard to verify that $F_0(l)G_0(l)^{-2}$ is non-decreasing when $0<p\leq \frac{n}{3n-4}$. By definitions of $G_1$ and $F_1$, when $0< p\leq \frac{n}{3n-3}$, one has
\begin{equation*}
	G_1''(l)=\frac{d}{dl}\frac{n-1}{\tan(I_1^{-1}(l))}\leq 0, \quad
	G_1(l)-3plG_1'(l)\geq 0.
\end{equation*}
Hence we have
\begin{equation*}
	\begin{aligned}
		\frac{d}{dl}\left(F_1'(l)G_1(l)-2F_1(l)G_1'(l)\right)&=F_1''(l)G_1(l)-F_1'(l)G_1'(l)-2F_1(l)G_1''(l)\\
		&\geq F_1''(l)G_1(l)-F_1'(l)G_1'(l)\\
		&=\frac{1}{p}l^{\frac{1}{p}-1}G_1(l)^{-2}\left(G_1(l)-3plG_1'(l)\right)\int_{0}^{l}{f^*(t) dt}+l^{\frac{1}{p}}G_1(l)^{-1}f^*(l)\\
		&\geq 0.
	\end{aligned}	
\end{equation*}
Noticing that $F_1'(0)G_1(0)-2F_1(0)G_1'(0)=0$, we have
\begin{equation*}
	F_1'(l)G_1(l)-2F_1(l)G_1'(l)\geq 0.
\end{equation*}
Thus, $F_1(l)G_1(l)^{-2}$ is non-decreasing since
\begin{equation*}
	(F_1G_1^{-2})'(l)=G_1(l)^{-3}(F_1'(l)G_1(l)-2F_1(l)G_1'(l))\geq 0.
\end{equation*}
When $n=2$, $0<p\leq 1$, we have
\begin{equation*}
	(G_1G_1')'(l)=2\pi\operatorname{\alpha_1}\frac{d \cos(I_1^{-1}(l))}{dl}\leq 0,\quad
	G_1(l)-2plG_1'(l)\geq 0.
\end{equation*}
Hence we deduce
\begin{equation*}
	\begin{aligned}
		\frac{d}{dl}\left(F_1'(l)G_1(l)^2-2F_1(l)G_1(l)G_1'(l)\right)&=F_1''(l)G_1(l)^2-2F_1(l)(G_1G_1')'(l)\\
		&\geq F_1''(l)G_1(l)^2\\
		&=\frac{1}{p}l^{\frac{1}{p}-1}G_1(l)^{-1}\left(G_1(l)-2plG_1'(l)\right)\int_{0}^{l}{f^*(t) dt}+l^{\frac{1}{p}}f^*(l)\\
		&\geq 0.
	\end{aligned}	
\end{equation*}
Together with $F_1'(0)G_1(0)^2-2F_1(0)G_1(0)G_1'(0)=0$, we have
\begin{equation*}
	F_1'(l)G_1(l)^2-2F_1(l)G_1(l)G_1'(l)\geq 0.
\end{equation*}
Then $F_1(l)G_1(l)^{-2}$ is non-decreasing by
\begin{equation*}
	(F_1G_1^{-2})'(l)=G_1(l)^{-4}(F_1'(l)G_1(l)^2-2F_1(l)G_1(l)G_1'(l))\geq 0.
\end{equation*}
Since $F_\kappa(l)G_\kappa(l)^{-2}$ is non-decreasing, by \eqref{lem2-1} we have
\begin{equation}
	\int_{0}^{\tau}{tF_\kappa(\mu(t)) dt}\leq -\int_{0}^{\tau}{t dH_\kappa(\mu(t))}+\frac{F_\kappa(|\Omega|_{g}) G_\kappa(|\Omega|_{g})^{-2}}{2\beta^2}\left(\int_{0}^{|\Omega|_{g}}{f^*(s) ds}\right)^2,
\end{equation}
where $H_\kappa(l)=\int_{0}^{l}{F_\kappa(w)G_\kappa(w)^{-2}\int_{0}^{w}{f^*(s) ds} dw}$. Integrating by parts and applying \eqref{Gronwall's inequality},  one gets
\begin{equation}\label{4.15}
	\begin{aligned}
		\int_{0}^{\tau}{F_\kappa(\mu(t)) dt}+H_\kappa(\mu(\tau))\leq\frac{1}{v_m}&\Bigg[\int_{0}^{v_m}{dt \int_{0}^{t}{F_\kappa(\mu(r)) dr} }+\int_{0}^{v_m}{H_\kappa(\mu(t)) dt}\\&+\frac{F_\kappa(|\Omega|_{g}) G_\kappa(|\Omega|_{g})^{-2}}{2\beta^2}\left(\int_{0}^{|\Omega|_{g}}{f^*(s) ds}\right)^2\Bigg].
	\end{aligned}
\end{equation}
Analogously, by \eqref{lem1-2} and \eqref{lem2-2}, we obtain
\begin{equation}
	\begin{aligned}
		\int_{0}^{\tau}{tF_\kappa(\operatorname{\alpha_\kappa}\phi(t)) dt}=& -\int_{0}^{\tau}{t\left(\operatorname{\alpha_\kappa}^{-1}\phi'(t) F_\kappa(\operatorname{\alpha_\kappa}\phi(t)) \tilde{G}_\kappa(\phi(t))^{-2} \int_{0}^{\operatorname{\alpha_\kappa}\phi(t)}{f^*(s) ds} \right) dt}\\
		&+\operatorname{\alpha_\kappa}^{-1}\int_{0}^{\tau}{\left[\frac{t}{\beta}\int_{\partial V_t\cap \partial \Omega^\sharp}{\frac{1}{v}d\mu_{g_\kappa}}\right] F_\kappa(\operatorname{\alpha_\kappa}\phi(t)) \tilde{G}_\kappa(\phi(t))^{-2} \int_{0}^{\operatorname{\alpha_\kappa}\phi(t)}{f^*(s) ds} dt}\\
		=&-\int_{0}^{\tau}{tdH_\kappa(\operatorname{\alpha_\kappa} \phi(t))}+\frac{F_\kappa(|\Omega|_{g}) G_\kappa(|\Omega|_{g})^{-2}}{2\beta^2}\left(\int_{0}^{|\Omega|_{g}}{f^*(s) ds}\right)^2.
	\end{aligned}	
\end{equation}
Thus,
\begin{equation}\label{4.17}
	\begin{aligned}
		\int_{0}^{\tau}{F_\kappa(\operatorname{\alpha_\kappa}\phi(t)) dt}+H_\kappa(\operatorname{\alpha_\kappa}\phi(\tau))=\frac{1}{v_m}&\Bigg[\int_{0}^{v_m}{dt \int_{0}^{t}{F_\kappa(\operatorname{\alpha_\kappa}\phi(r)) dr} }+\int_{0}^{v_m}{H_\kappa(\operatorname{\alpha_\kappa}\phi(t)) dt}\\
		&+\frac{F_\kappa(|\Omega|_{g}) G_\kappa(|\Omega|_{g})^{-2}}{2\beta^2}\left(\int_{0}^{|\Omega|_{g}}{f^*(s) ds}\right)^2\Bigg].
	\end{aligned}
\end{equation}
By \eqref{mf}, \eqref{4.15} and \eqref{4.17}, we come up with
\begin{equation}
	\int_{0}^{\tau}{F_\kappa(\mu(t)) dt}+H_\kappa(\mu(\tau))\leq \int_{0}^{\tau}{F_\kappa(\operatorname{\alpha_\kappa}\phi(t)) dt}+H_\kappa(\operatorname{\alpha_\kappa}\phi(\tau)).
\end{equation}
Letting $\tau\rightarrow +\infty$ and using \eqref{4.10}, \eqref{4.11}, we prove \eqref{thm:main1-2}, \eqref{thm:main1-4} and \eqref{thm:main1-5}. 
\end{proof}

\begin{proof}[Proof of Theorem \ref{thm:main2}]
	If $f\equiv 1$, then \eqref{52} becomes
	\begin{equation}
		\begin{aligned}
			\int_{0}^{\tau}{t\mu(t)^{\frac{1}{p}} dt}\leq& \int_{0}^{\tau}{-t\mu'(t) \mu(t)^{\frac{1}{p}+1} G_\kappa(\mu(t))^{-2} dt}\\
			&+\int_{0}^{\tau}{\left[\frac{t}{\beta}\int_{\partial U_t^e}{\frac{1}{u}d\mu_{g}}\right] \mu(t)^{\frac{1}{p}+1} G_\kappa(\mu(t))^{-2} dt}.
		\end{aligned}
	\end{equation}
It's easy to show that $l^{\frac{1}{p}+1}G_\kappa(l)^{-2}$ is non-decreasing when $\kappa=0,1$ and $0<p\leq \frac{n}{n-2}$, hence
\begin{equation}\label{74}
	\int_{0}^{\tau}{t\mu(t)^{\frac{1}{p}} dt}\leq
	-\int_{0}^{\tau}{tdF_\kappa(\mu(t))}+\frac{|\Omega|_{g}^{\frac{1}{p}+2} G_\kappa(|\Omega|_{g})^{-2}}{2\beta^2}.
\end{equation}	
Following the similar arguments in proof of \eqref{thm:main1-1}, we are able to prove \eqref{thm:main2-1}. 

Noticing that $lF_0(l)G_0(l)^{-2}$ is still non-decreasing when $0<p
\leq \frac{n}{n-2}$, we adopt same  as the proof of \eqref{thm:main1-2}, \eqref{thm:main1-4} and \eqref{thm:main1-5} and finish the proof of \eqref{thm:main2-2}.

Finally, we give the proof of the pointwise comparison result \eqref{est-point}.
For $n=2$,  $\kappa=0$ and $f\equiv 1$, \eqref{lem1-1} becomes
\begin{equation}\label{new-1}
	4\pi\akz\leq -\mu'(t)+\frac{1}{\beta}\int_{\partial U_t^e}{\frac{1}{u}d\mu_{g}},
\end{equation}
while \eqref{lem1-2} yields
\begin{equation}\label{new-2}
	4\pi=-\phi'(t)+\frac{1}{\beta}\int_{\partial V_t\cap\partial\Omega^\sharp}{\frac{1}{v}d\mu_{g_\kappa}}.
\end{equation}
Multiplying \eqref{new-1} by $t$, integrating form 0 to $\tau>v_m$ and applying \eqref{lem2-1}, we obtain
\begin{equation}
	2\pi\akz\tau^2\leq \int_{0}^{\tau}{-t\mu' dt}+\frac{\abs{\Omega}_g}{2\beta^2}.
\end{equation}
Meanwhile, by \eqref{new-2}, one has
\begin{equation}
	2\pi\akz\tau^2=\int_{0}^{\tau}{-t(\akz\phi(t))' dt}+\frac{\abs{\Omega}_g}{2\beta^2}.
\end{equation}
Thus,
\begin{equation}
	\int_{0}^{\tau}{-t\mu' dt}\geq\int_{0}^{\tau}{-t(\akz\phi(t))' dt}.
\end{equation}
Integrating by parts and using \eqref{Gronwall's inequality}, for all $\tau>v_m$, we have
\begin{equation}
	\mu(\tau)-\akz\phi(\tau)\leq\frac{1}{v_m}\int_{0}^{v_m}{\mu(t)-\akz\phi(t) dt}.
\end{equation}
Together with \eqref{mf}, we find out that for all $\tau\geq 0$, 
\begin{equation}
	\mu(\tau)\leq\akz\phi(\tau).
\end{equation}
By definition of Schwarz rearrangement, we complete the proof of \eqref{est-point}.
	\end{proof}

  \section{Proof of Theorem \ref{thm:eig}}\label{pf:BD}
  
  For $n=2$, it is independently interesting to prove Bossel-Daners inequality by using our Theorem \ref{thm:main1}. For $n\geq 2$, the methods in \cite{Bossel} and \cite{D1}  can be adopted  to prove the Bossel-Daners inequality. For bounded domains in Riemannian manifolds with $\Ric(g)\geq (n-1)$, the inequality \eqref{thm:eig-1} has already been proven in \cite{CCL2021}. Here we only give the proof of \eqref{thm:eig-1} for complete noncompact Riemannian manifolds with nonnegative Ricci curvature and positive asymptotic volume growth.
  
  \begin{proof}[Proof of Theorem \ref{thm:eig}] Following the idea in \cite{Kesavan88}, we can prove the  Bossel-Daners inequality for the first eigenvalue of  Robin Laplacian in dimension $2$. 	Let $u$ and $z$ be solutions to \eqref{eig-1} and \eqref{eig-2} respectively.
When $n=2$, by choosing $p=1$ in \eqref{thm:main1-2} and \eqref{thm:main1-5}, we have
\begin{equation}
	\int_{\Omega^\sharp}{(u^\sharp)^2 dV_{g_\kappa}}=\ak^{-1}\int_{\Omega}{u^2 dV_{g}}\leq \int_{\Omega^\sharp}{z^2 dV_{g_\kappa}}.
\end{equation}
By Cauchy-Schwarz inequality, one has
\begin{equation}
	\int_{\Omega^\sharp}{u^\sharp z dV_{g_\kappa}}\leq \int_{\Omega^\sharp}{z^2 dV_{g_\kappa}}.
\end{equation}
By Rayleigh quotient for eigenvalues, we obtain
\begin{equation}
	\begin{aligned}
		\lambda_{1,\beta}(\Omega)&=\frac{\int_{\Omega^\sharp}{|\nabla_{g_\kappa} z|_{g_\kappa}^2 dV_{g_\kappa}}+\beta\int_{\partial \Omega^\sharp}{z^2 d\mu_{g_\kappa}}}{\int_{\Omega^\sharp}{u^\sharp z dV_{g_\kappa}}}\\
		&\geq \frac{\int_{\Omega^\sharp}{|\nabla_{g_\kappa} z|_{g_\kappa}^2 dV_{g_\kappa}}+\beta\int_{\partial \Omega^\sharp}{z^2 d\mu_{g_\kappa}}}{\int_{\Omega^\sharp}{z^2 dV_{g_\kappa}}}\geq \lambda_{1,\beta}(\Omega^\sharp).
	\end{aligned}	
\end{equation}	
Now we are in position to give  the  Bossel-Daners inequality for  $n\geq2$ and $\kappa=0$.
Let $B^0_R$ be the geodesic ball of radius $R$ in Euclidean space $(\R(0),g_0)$ such that $\abs{\Omega}_g=\operatorname{\alpha_0}\abs{B^0_R}_{g_0}$.  Let $u_0$ be the eigenfunction associated to  $\lambda_{1,\beta}(B^0_R)$  on geodesic ball $B^0_R$, i.e. 
\begin{equation}\label{robin-ball}
	\begin{cases}
		u_0^{\prime\prime}+\frac{n-1}{r}u_0^\prime+\lambda_{1,\beta}(B^0_R)u_0=0, \quad \text{on}~B^0_R,\\
		u_0^{\prime}(0)=0, \quad u_0^{\prime}(R)+\beta u_0(R)=0.
	\end{cases}	
\end{equation}
It is known that the function $u_0(r)$ is positive for $r\in [0,R]$.
Define  $v(r)=(\ln u_0)^\prime(r)$. According to Lemma 2.1 and Proposition 2.2  in \cite{CCL2021} (See also \cite{D1,BD}), the function $v(r)=\left(\ln u_0\right)^\prime(r)$ is strictly decreasing  and then $0\leq -v(r)<\beta$ for $r\in [0,R)$ and the first eigenvalue $\lambda_{1,\beta}(B^0_R)$ of Robin Laplacian  is strictly decreasing  with respect to radius $R$.
Let  $u$ be the first eigenfunction of Laplacian equation on $\Omega$. We can choose $u>0$ and normalize it such that
\begin{equation*}
	\begin{aligned}
		&\underset{x\in \overline{\Omega}}{\max} \, u(x)=1,\quad
		\underset{x\in \overline{\Omega}}{\min} \, u(x)=m.
	\end{aligned}
\end{equation*}
For $t \in (m,1) $,  set
\begin{equation}\label{set}
		U_t =\{x \in \Omega: \  u(x) >t\},\quad 
		\partial U_t^i=\{x\in\Omega: u(x)=t\},\quad
		\partial U_t^e=\partial U_t\cap \partial \Omega=\{x\in\partial\Omega: u(x)\geq t \}.
\end{equation}
Define the admissible subset $M _ { \beta } (\Omega)$ of $C(\Omega)$ by
\begin{equation}\label{Adm-M}
	M _ { \beta }(\Omega)  = \left\{ \varphi\in C ( \Omega ) : \,\varphi(x)\geq 0, \limsup _ { x \rightarrow z } \varphi ( x ) \leq \beta, \quad \text { for all } z \in \partial \Omega \right\}.
\end{equation}
For $\varphi\in M_\beta(\Omega)$, we can define the functional  $H_{\Omega}(U_t,\varphi)$ as in \cite{Bossel} and \cite{D1}
\begin{equation}\label{DB-functional}
	H_{\Omega}(U_t,\varphi)=\frac{1}{\abs{U_t}_g}\left(\beta \abs{\partial U_t^e}_g+\int_{\partial U_t^i}\varphi d\mu_g-\int_{U_t}\varphi^2dV_g\right).
\end{equation}
For the first eigenfunction $u$, one can infer that $\frac{\abs{\n u}_g}{u}\in  M _ { \beta }(\Omega)$ and 
\begin{equation}\label{DB-functional1}
	\ld=H_{\Omega}\left(U_t,\frac{\abs{\n u}_g}{u}\right),\qquad \text{for almost all}\, t\in (m,1).
\end{equation}
For $\varphi\in M_\beta(\Omega)$, set
\begin{equation*}
	w=\varphi-\frac{\abs{\n u}_g}{u},\qquad \text{and}\qquad F(t)=\int_{U_t}w\frac{\abs{\n u}_g}{u}dV_g.
\end{equation*}
From Lemma 3.1 in \cite{CCL2021} (See also \cite{BD,D1}), the functional $H_{\Omega}$ satisfies
\begin{equation}\label{equ3.0}
	H_{\Omega}(U_t,\varphi)=\ld-\frac{1}{\abs{U_t}_g}\left(\frac{1}{t}\frac{d}{dt}\left(t^2F(t)\right)+\int_{U_t}w^2dV_g\right)
\end{equation}
for almost all $t\in (m,1)$. Furthermore, suppose that $\varphi\neq\frac{\abs{\n u}_g}{u} $, then there exits  a set $S\subset (m,1)$ with $\abs{S}>0$ such that
\begin{equation}\label{DB-functional2}
	\ld>H_{\Omega}(U_t,\varphi),\qquad \text{for all}\quad t\in S.
\end{equation}
For  $t\in (m,1]$, we define the ball $B^0_{r(t)}$ with radius $r(t)$  such that $\abs{U_t}_g=\operatorname{\alpha_0}\abs{B^0_{r(t)}}_{g_0}$. For $x \in \partial U_t^i$ and $t\in (m,1]$, we define
\begin{equation}\label{funct-test}
	\varphi(x)=-v(r(t))=\frac{-u_0^\prime(r(t))}{u_0(r(t))}.
\end{equation}
It is easy to check that  $\varphi:\Omega \rightarrow (0, \infty)$   is a measurable function and  
\begin{equation}\label{}
	U_t=\{x \in \Omega: \  u(x) >t\}=\{x \in \Omega: \  \varphi(x)<-v(r(t))\}
\end{equation}
is open in $\Omega$ for $t\in (m,1)$.  From  the construction of $\varphi$ in \eqref{funct-test} and the monotonic decreasing property of $v(r)=\left(\ln u_0\right)^\prime(r)$, we notice that
\begin{equation*}
	0\leq \varphi(x)<\beta, \qquad\text{for}\quad  x\in \Omega,
\end{equation*}
which implies that $\varphi\in M_{\beta}(\Omega)$. Since $\abs{U_t}_g=\operatorname{\alpha_0}\abs{B^0_{r(t)}}_{g_0}$,
we have
\begin{equation}\label{Cavalieri}
	\int_{U_t}\varphi^2dV_g=\operatorname{\alpha_0} \int_{B^0_{r(t)}}(-v(r(t))^2dV_{g_0}, \quad 
	\quad \text{for}\, t\in (m,1].
\end{equation}
Moreover, by isoperimetric inequality \eqref{isoperi} and $\varphi=-v(r(t))\leq \beta$, we obtain that
\begin{equation}\label{iso1}
	\begin{aligned}
		\operatorname{\alpha_0}\abs{\partial B^0_{r(t)}}_{g_0}(-v(r(t)))\leq& (-v(r(t)))\abs{\partial U_t}_g\\
		\leq &\int_{\partial U_t^i}\varphi(x)d\mu_g+\beta\abs{\partial U_t^e}_g, \qquad \text{for} \quad t\in (m,1),
	\end{aligned}
\end{equation}
From \eqref{funct-test}  and \eqref{iso1}, we conclude
\begin{equation*}
	\begin{aligned}
		\ldb&=H_{B^0_R}(B^0_{r(t)},-v(r(t)))\\
		&=\frac{1}{\operatorname{\alpha_0}\abs{B^0_{r(t)}}_{g_0}}\left(\operatorname{\alpha_0}\int_{\partial B^0_{r(t)}} (-v(r(t)))d\mu_{g_0}-\operatorname{\alpha_0}\int_{B^0_{r(t)}} (-v(r(t)))^2dV_{g_0}\right)\\
		&\leq  \frac{1}{\abs{U_t}_g}\left(\beta\abs{\partial U_t^e}_g+\int_{\partial U_t^i}\varphi(x)d\mu_g-\int_{U_t}\varphi^2dV_g\right)\\
		&\leq  \ld,
	\end{aligned}
\end{equation*}
which completes the proof of the inequality \eqref{thm:eig-1}. When equality occurs in \eqref{thm:eig-1}, the above inequalities become equalities. Therefore, according to \cite{Brendle}, the equality holds in \eqref{thm:eig-1} if and only if $(M,g)$ is isometric to $(\R(0),g_0)$ and $\Omega$ is isometric to a ball $B^0_R$, where $g_0$ is the canonical metric of Euclidean space.

\end{proof}

\providecommand{\bysame}{\leavevmode\hbox
	to3em{\hrulefill}\thinspace}

\vspace{1cm}

\begin{flushleft}
	Daguang Chen,
	E-mail: dgchen@tsinghua.edu.cn\\
	Haizhong Li,
	E-mail:	lihz@tsinghua.edu.cn\\
	Yilun  Wei,
	E-mail:	weiyl19@mails.tsinghua.edu.cn\\
	Department of Mathematical Sciences, Tsinghua University, Beijing, 100084, P.R. China 	
	
\end{flushleft}

\end{document}